\documentclass{gtart_h}

\def\ifplaintex{\expandafter\ifx\csname documentclass\endcsname\relax}

\def\gtp{{\mathsurround=0pt\it $\cal G\mskip-2mu$eometry \&\ 
$\cal T\!\!$opology $\cal P\!$ublications}}  

\def\recd{{\small Received:\qua\receiveddate\ifx\reviseddate\relax
\else\qquad Revised:\qua\reviseddate\fi\par}} 

\def\lognumber#1{\def\thelognumber{#1}}
\def\volumenumber#1{\def\thevolumenumber{#1}}
\def\volumeyear#1{\def\thevolumeyear{#1}}
\def\papernumber#1{\def\thepapernumber{#1}}
\def\pagenumbers#1#2{\def\startpage{#1}\def\finishpage{#2}}
\def\published#1{\def\publishdate{#1}}

\def\received#1{\def\receiveddate{#1}}

\def\accepted#1{\def\accepteddate{#1}}

\long\def\asciiabstract#1{\long\def\theasciiabstract{#1}}

\let\\\par\let\thelognumber\relax\let\thevolumenumber\relax
\let\thepapernumber\relax\let\thevolumeyear\relax\let\startpage\relax
\let\finishpage\relax\let\publishdate\relax\let\receiveddate\relax
\let\reviseddate\relax\let\accepteddate\relax\let\theasciititle\relax
\let\theasciiauthors\relax
\let\theasciiabstract\relax

\let\theasciiemail\relax

\ifplaintex
\font\logobig=cmssbx10 scaled 3836
\font\logomed=cmssbx10 scaled 2557
\else
\font\logobig=cmssbx10 scaled 4200
\font\logomed=cmssbx10 scaled 2800
\fi

\long\def\makeagttitle{   
\count0=\startpage
\agt\hfill      
\hbox to 45truept{\vbox to 0pt{\vglue -13truept{\logomed A\kern -.37em{\logobig 
T}\kern -.38em G}\vss}\hss}
\break
{\small Volume \thevolumenumber\ (\thevolumeyear)
\startpage--\finishpage\nl
Published: \publishdate}

\vglue .25truein

{\parskip=0pt\leftskip 0pt plus
1fil\def\\{\par\smallskip}{\Large\bf\thetitle}\par\medskip} \vglue
0.05truein

{\parskip=0pt\leftskip 0pt plus 1fil\def\\{\par}{\sc\theauthors}
\par\medskip}%
 
\vglue 0.03truein 

{\small\leftskip 25truept\rightskip 25truept{\bf Abstract}\stdspace\theabstract

{\bf AMS Classification}\stdspace\theprimaryclass
\ifx\thesecondaryclass\relax\else; \thesecondaryclass\fi\par
{\bf Keywords}\stdspace \thekeywords\par}\vglue 7truept

}   

\ifplaintex
\font\phead=cmsl9 scaled 950
\font\pnum=cmbx10 scaled 913
\font\pfoot=cmsl9 scaled 950
\headline{\vbox to 0pt{\vskip -4.5mm\line{\small\phead\ifnum
\count0=\startpage ISSN 1472-2739 (on-line) 1472-2747 (printed)
\hfill {\pnum\folio}\else\ifodd\count0\def\\{ }%
\ifx\theshorttitle\relax\thetitle\else\theshorttitle\fi\hfill{\pnum\folio}
\else\def\\{ and }{\pnum\folio}\hfill\ifx\theshortauthors\relax\theauthors
\else\theshortauthors\fi\fi\fi}\vss}}
\footline{\vbox to 0pt{\vglue 0mm\line{\small\pfoot\ifnum\count0=\startpage
\copyright\ \gtp\hfill\else
\agt, Volume \thevolumenumber\ (\thevolumeyear)\hfill\fi}\vss}}
\else
\font=cmsl9 scaled 1050
\font\lnum=cmbx10 
\font=cmsl9 scaled 1050
\makeatletter
\def\@oddhead{{\small\lhead\ifnum\count0=\startpage ISSN 1472-2739 
(on-line) 1472-2747 (printed)\hfill {\lnum\number\count0}\else\ifodd\count0
\def\\{ }\ifx\theshorttitle\relax \thetitle \else\theshorttitle\fi\hfill
{\lnum\number\count0}\else\def\\{ and }{\lnum\number\count0}
\hfill\ifx\theshortauthors\relax 
\theauthors\else\theshortauthors\fi\fi\fi}}\def\@evenhead{\@oddhead}
\def\@oddfoot{\small\lfoot\ifnum\count0=\startpage\copyright\ \gtp\hfill\else
\agt, Volume \thevolumenumber\ (\thevolumeyear)\hfill\fi}
\def\@evenfoot{\@oddfoot}
\makeatother
\fi
\let\maketitlepage\makeagttitle
\let\makeshorttitle\maketitlepage
\let\maketitle\maketitlepage

\newwrite\gtoutfile
\long\gdef\makeheadfile{  
{\def\\{, }\def\s{ }
\immediate\openout\gtoutfile head.xxx
\immediate\write\gtoutfile{Proxy-for: \ifx\theasciiauthors\relax
\theauthors\else\theasciiauthors\fi\s<\ifx\theasciiemail\relax\theemail\else\theasciiemail\fi>}
\immediate\write\gtoutfile{\noexpand\\}
\immediate\write\gtoutfile{Authors: \ifx\theasciiauthors\relax
\theauthors\else\theasciiauthors\fi}
{\def\\{ }\immediate\write\gtoutfile{Title: \ifx\theasciititle\relax
\thetitle\else\theasciititle\fi}}
\immediate\write\gtoutfile{Subj-class: GT or SG, GR etc}
\immediate\write\gtoutfile{MSC-class: \theprimaryclass\ifx\thesecondaryclass\relax\else, \thesecondaryclass\fi}
\immediate\write\gtoutfile{Journal-ref: Algebr. Geom. Topol. \thevolumenumber\s
(\thevolumeyear) \startpage-\finishpage}
\immediate\write\gtoutfile{Comments: Published by Algebraic and
Geometric Topology at}
\immediate\write\gtoutfile{\s\s\s  http://www.maths.warwick.ac.uk/agt/AGTVol\thevolumenumber/agt-\thevolumenumber-\thepapernumber.abs.html}
\immediate\write\gtoutfile{\noexpand\\}
\immediate\write\gtoutfile{}
\ifx\theasciiabstract\relax
\immediate\write\gtoutfile{\theabstract}\else
\immediate\write\gtoutfile{\theasciiabstract}\fi
\immediate\write\gtoutfile{}
\immediate\write\gtoutfile{\noexpand\\}
\immediate\write\gtoutfile{}
\immediate\closeout\gtoutfile}}  

\def\maketitlepage{\makeagttitle\makeheadfile}
\let\makeshorttitle\maketitlepage
\let\maketitle\maketitlepage

\lognumber{15}
\volumenumber{4}
\volumeyear{2004}
\papernumber{15}
\published{25 April 2004}
\pagenumbers{273}{296}
\received{28 March 2003}
\accepted{9 February 2004}

\usepackage{graphicx, amsmath, amssymb, amscd, psfrag}

\theoremstyle{definition}
\newtheorem{definition}{Definition}[section]
\newtheorem{example}[definition]{Example}
\newtheorem{lemma}[definition]{Lemma}
\theoremstyle{plain}
\newtheorem{proposition}[definition]{Proposition}
\newtheorem{theorem}[definition]{Theorem}
\newtheorem{corollary}[definition]{Corollary}
\newtheorem*{main}{Main Theorem}
\newtheorem*{maint}{Main Theorem}

\newcommand{\rw}{\rightarrow}
\newcommand{\lw}{\leftarrow}
\newcommand{\zt}{\mathbb{Z}_2}
\newcommand{\zz}{\mathbb{Z}}
\newcommand{\nn}{\mathbb{N}}
\newcommand{\rr}{\mathbb{R}}

\begin{document}

\newcounter{bean}

\title{A lower bound to the action\\dimension of a group}
\author{Sung Yil Yoon}
\address{110 8th Street RPI, Troy, NY 12180, USA}
\email{yoons@rpi.edu}
\primaryclass{20F65} 
\secondaryclass{57M60}

\keywords{Fundamental group, contractible manifold, action dimension, 
embedding obstruction}

\begin{abstract}
The \emph{action dimension} of a discrete group
$\Gamma$, $actdim(\Gamma)$, is defined to be the smallest integer
$m$ such that $\Gamma$ admits a properly discontinuous action on
a contractible $m$--manifold. If no such $m$ exists, we define
$actdim(\Gamma) \equiv \infty$.
Bestvina, Kapovich, and Kleiner used Van Kampen's theory of
embedding obstruction to provide a lower bound to the action
dimension of a group. In this article, another lower bound to the
action dimension of a group is obtained by extending their work,
and the action dimensions of the fundamental groups of certain
manifolds are found by computing this new lower bound.
\end{abstract}
\asciiabstract{%
The action dimension of a discrete group G, actdim(G), is defined to
be the smallest integer m such that G admits a properly discontinuous
action on a contractible m-manifold. If no such m exists, we define
actdim(G) = infty.  Bestvina, Kapovich, and Kleiner used Van Kampen's
theory of embedding obstruction to provide a lower bound to the action
dimension of a group. In this article, another lower bound to the
action dimension of a group is obtained by extending their work, and
the action dimensions of the fundamental groups of certain manifolds
are found by computing this new lower bound.}

\maketitle

\section{Introduction}

Van Kampen constructed an $m$--complex that cannot be embedded into
$\rr^{2m}$~\cite{vK}.  A more modern approach to Van Kampen's
theory of embedding obstruction uses co/homology theory.
To see the main idea of this co/homology theoretic approach, let
$K$ be a simplicial complex and $|K|$ denote its geometric
realization. Define  the \emph{deleted product}
\[ \tilde{|K|} \equiv \{ (x,y) \in |K| \! \times \! |K|\  |\  x \neq y \} \]
such that $\zt$ acts on $\tilde{|K|}$ by exchanging factors.
Observe that there exists a two-fold covering $\tilde{|K|} \rw \tilde{|K|}/\zt$
with the following classifying map:
\begin{equation*}
\begin{CD}
\tilde{|K|}  @>\tilde{\phi}>>    S^{\infty} \\
     @VVV                           @VVV \\
\tilde{|K|}/\zt @>\phi>>     \rr P^{\infty}
\end{CD}
\end{equation*}
Now let $\omega^m \in H^m(\rr P^{\infty};\zt)$ be the nonzero class.
If $\phi^*(w^m) \neq 0$ then $|K|$ cannot be embedded into $\rr^m$.
That is, there is $\Sigma^m \in H_m(\tilde{|K|}/\zt; \zt)$ such that
$\langle \phi^*(w^m), \Sigma \rangle \neq 0$.

A similar idea was used to obtain a lower bound to the
\emph{action dimension} of a discrete group $\Gamma$~\cite{BKK}.
Specifically, the \emph{obstructor dimension}
of a discrete group $\Gamma$, $obdim(\Gamma)$,  was defined by
considering an $m$--\emph{obstructor} $K$ and a \emph{proper},
Lipschitz, \emph{expanding} map
\[ f\colon\thinspace cone(K)^{(0)} \rightarrow \Gamma .\]
And it was shown that
\[ obdim(\Gamma) \leq actdim(\Gamma). \]
See~\cite{BKK} for details.
An advantage of considering $obdim(\Gamma)$ becomes clear when
$\Gamma$ has well-defined boundary $\partial \Gamma$, for example,
when $\Gamma$ is $CAT(0)$ or torsion free hyperbolic.  In these cases,
if an $m$--obstructor $K$ is contained in
$\partial \Gamma$ then $m+2 \leq obdim(\Gamma)$.

If $\Gamma$ acts on a contractible $m$--manifold $W$ properly discontinuously
and cocompactly, then it is easy to see that $actdim(\Gamma)=m$.
For example, let $M$ be a Davis manifold.
That is, $M$ is a closed, aspherical,
four-dimensional manifold whose universal cover $\tilde{M}$ is not
homeomorphic to $\rr^4$.  We know that $actdim (\pi_1(M) )=4$. However, it
is not easy to see that $obdim (\pi_1(M)) =4$.   The goal of this article
is to generalize the definitions of obstructor and obstructor dimension.
To do so, we define \emph{proper obstructor} (Definition~\ref{pobst}) and
\emph{proper obstructor dimension} (Definition~\ref{pobdim}.)
The main result is the following.
\begin{main}
The proper obstructor dimension of $\Gamma$ $\leq actdim(\Gamma)$.
\end{main}

\indent
As applications we will answer the following problems:

\begin{itemize}
\item Suppose $W$ is a closed aspherical manifold and $\tilde{W}$ is its
universal cover so that $\pi_1(W)$ acts on $\tilde{W}$ properly
discontinuously and cocompactly.  We show that $\tilde{W}$ in this
case is indeed an $m$--proper obstructor and $pobdim(\pi_1(W))=m$.

\item Suppose $W_i$ is a compact aspherical $m_i$--manifold with all boundary
components aspherical and incompressible, $i=1, ..., d$. (Recall that a
boundary component $N$ of a manifold $W$ is called \emph{incompressible}
if $i_*\colon\thinspace\pi_j(N) \rw \pi_j(W)$ is injective for $j \geq 1$.)   Also assume
that for each $i$, $1\leq i \leq d$, there is a component of $\partial W_i$,
call it $N_i$, so that $|\pi_1(W_i)\! : \! \pi_1(N_i)|>2$.
Let $G=\pi_1(W_1) \! \times ... \times \!  \pi_1(W_d)$.   Then
\[actdim(G)=m_1+ ... + m_d .\]

\end{itemize}

The organization of this article is as follows.
In Section~\ref{po}, we define proper obstructor.
The \emph{coarse Alexander duality} theorem by Kapovich and
Kleiner~\cite{KK}, is used to construct the first main example of
proper obstructor in Section~\ref{ca}. Several examples of
proper obstructors are constructed in Section~\ref{pod1}.
Finally, the main theorem is proved and the above problems are
considered in Sections~\ref{pod2}.

\section{Proper obstructor}
\label{po}

To work in the
PL--category we define \emph{simplicial deleted product}
\[\tilde{K} \equiv \{ \sigma \! \times \! \tau \in K\! \times \! K\  |\
\sigma \cap \tau =\emptyset \} \] such that  $\zt$  acts on $\tilde{K}$
by exchanging factors.  It is known that $\tilde{|K|}/\zt$ $(\tilde{|K|})$
is a deformation retract of $\tilde{K}/\zt$ $(\tilde{K})$,
see~\cite[Lemma 2.1]{Sh}.
Therefore, WLOG, we can use $H_m(\tilde{K}/\zt ; \zt)$
instead of $H_m(\tilde{|K|}/\zt ; \zt)$.

\textsl{Throughout the paper, all homology groups are taken with
$\zt$--coefficients unless specified otherwise}.

To define proper obstructor, we need to consider several definitions
and preliminary facts.
\begin{definition}
A proper map $h \colon\thinspace A \rw B$ between proper metric spaces is
\emph{uniformly proper} if there is a proper function
$\phi \colon\thinspace [0,\infty) \rw [0,\infty)$ such that
\[ d_B(h(x), h(y) ) \geq \phi( d_A(x,y) ) \]
for all $x,y \in A$. (Recall that a metric space is said to be
\emph{proper} if any closed metric ball is compact, and a map is
said to be \emph{proper} if the preimages of compact sets are
compact.)
\end{definition}
\noindent
Let $W$ be a contractible $m$--manifold and define
\[W_0 \equiv \{(x,y) \in W \! \times \! W\  |\  x \neq y \}. \]
Consider a uniformly proper map $\beta \colon\thinspace Y \rw W$ where $Y$ is a
simplicial complex and $W$ is a contractible manifold.  Since
$\beta$ is uniformly proper, we can choose $r>0$ such that
$\beta(a) \neq \beta(b)$ if $d_Y(a,b) >r$.  Note that
$\beta$ induces an equivariant map:
\[\bar{\beta} \colon\thinspace \{ (y,y') \in Y\! \times \! Y\ |\ d_Y (y,y') >r \} \rw W_0 \]
As we work in the PL--category we make the following definition.
\begin{definition}
\label{nbd}
If $K \subset Y$ is a subcomplex and $r$ is a positive integer then
we define the \emph{combinatorial $r$-tubular neighborhood of $K$},
denoted by $N_r(K)$, to be $r$-fold iterated closed star neighborhood of $K$.

\end{definition}
\noindent
Recall that when $Y$ is a simplicial complex, $|Y|\! \times \! |Y|$ can be
triangulated so that each cell $\sigma \! \times \! \tau$ is a subcomplex.
Let $d \colon\thinspace Y \rw Y^2$ be the diagonal map, $d(\sigma)=(\sigma, \sigma)$,
where $Y^2$ is triangulated so that $d(Y)$ is a subcomplex.
Define \[ Y_r \equiv Cls(\, Y^2- N_r(d(Y))\, ). \]
Note that a uniformly proper map $\beta \colon\thinspace Y \rw W^m$ induces an equivariant
map $\bar{\beta} \colon\thinspace Y_r \rw W_0 \simeq S^{m-1}$ for some $r>0$.

\begin{definition}[(Essential $\zt\!-\!m\!-\!$cycle)]
\label{ess-cycle}
An \emph{essential $\zt\!-\!m\!-\!$cycle} is a pair \linebreak
$(\tilde{\Sigma}^m, a)$ satisfying the following conditions:

\begin{itemize}
\item[(i)] $\tilde{\Sigma}^m$ is a finite simplicial complex such that
$|\tilde{\Sigma}^m|$ is a union of $m$--simplices and every $(m-1)$--simplex
is the face of an even number of $m$--simplices.
\item[(ii)] $a \colon\thinspace \tilde{\Sigma}^m \rw \tilde{\Sigma}^m$ is a free involution.
\item[(iii)] There is an equivariant map $\varphi \colon\thinspace  \tilde{\Sigma}^m \rw S^m$
with $deg(\varphi)=1 (mod \ 2)$.
\end{itemize}
\end{definition}

\indent
Some remarks are in order.
\begin{enumerate}

\item We recall how to find $deg(\varphi)$.  Choose a simplex $s$ of $S^m$ and
let $f$ be a simplicial approximation to $\varphi$.  Then
$deg(\varphi)$ is the number of $m$--simplices of $\tilde{\Sigma}^m$
mapped into $s$ by $f$.

\item  Let $\tilde{\sigma}$ be the sum of all $m$--simplices of
$\tilde{\Sigma}^m$.  Condition (i) of Definition~\ref{ess-cycle}
implies that $[\tilde{\sigma}] \in H_m(\tilde{\Sigma}^m)$.
We call $[\tilde{\sigma}]$ the fundamental class of $\tilde{\Sigma}^m$.

\item Let $\tilde{\Sigma}^m /\zt \equiv \Sigma^m$ and consider a two-fold
covering $q \colon\thinspace \tilde{\Sigma}^m \rw \Sigma^m$.  As $\varphi$ is equivariant
it induces $\bar{\varphi} \colon\thinspace \Sigma^m \rw \rr P^{m}$.   Let $deg_2(\varphi)$ denote
$deg(\varphi) (mod\  2)$.    

Note that
$deg_2(\varphi)=\langle \bar{\varphi}^{*} (w^m),  q\tilde{\sigma} \rangle$ where
$w^m \in H^m (\rr P^m; \zt)$ is the nonzero element.
If $\varphi \colon\thinspace \tilde{\Sigma}^m \rw S^m$ is an equivariant map then $deg_2(\varphi)=1$.
To see this, we prove the following proposition.
\end{enumerate}

\begin{proposition}
\label{well-prop}
Suppose a map $\varphi \colon\thinspace \tilde{\Sigma}^m \rw S^m$
is equivariant.  Then $$deg_2 (\varphi) =1.$$
\end{proposition}
\proof
Consider the classifying map and the commutative diagram for a
two-fold covering $q \colon\thinspace \tilde{\Sigma}^m  \rw \Sigma^m$:
\begin{equation*}
\begin{CD}
\tilde{\Sigma}^m   @>\phi >>    S^{\infty} \\
     @V q VV                       @V p VV \\
\Sigma^m           @>\bar{\phi}>> \rr P^{\infty}
\end{CD}
\end{equation*}
We also have:
\begin{equation*}
\begin{CD}
\tilde{\Sigma}^m   @>\varphi >>   S^m        @>i>>   S^{\infty} \\
     @V q VV                 @V p VV                  @V p VV      \\
\Sigma^m    @>\bar{\varphi} >> \rr P^m  @>i>>   \rr P^{\infty}
\end{CD}
\end{equation*}
Because $S^{\infty} \rw \rr P^{\infty}$ is the classifying covering,
$i \circ \varphi \simeq \phi$ and
$i \circ \bar{\varphi} \simeq \bar{\phi}$.
Observe that \[ deg_2(\varphi) =
\langle {\bar{\varphi}}^* w^m ,\  q\tilde{\sigma} \rangle =
\langle {(i \circ \bar{\varphi})}^*  w^m_{\infty},\  q \tilde{\sigma} \rangle \]
where $0 \neq w^m_{\infty} \in H^m (\rr P^{\infty} )$.
But, since $i \circ \bar{\varphi} \simeq \bar{\phi}$,
$$ \langle {(i \circ \bar{\varphi})}^*  w^m_{\infty},\ 
q \tilde{\sigma} \rangle = \langle \bar{\phi}^*  w^m_{\infty},\ 
q \tilde{\sigma} \rangle = deg_2(\phi).\eqno{\qed} $$

Now we modify the definition of obstructor.

\begin{definition}[(Proper obstructor)]
\label{pobst}
Let $T$ be a contractible\footnote{Contractibility is necessary for Proposition 5.2}
simplicial complex.  Recall that $T_r \equiv
Cls(\, T^2 - N_r(d(T))\,)$ where $N_r(d(T))$ denotes the
$r$--tubular neighborhood of the image of the diagonal map $d \colon\thinspace T \rw T^2$.
Let $m$ be the largest integer such that for any  $r>0$, there exists an essential
$\zt$--$m$--cycle $(\tilde{\Sigma}^m, a)$ and a
$\zt$--equivariant map $f \colon\thinspace \tilde{\Sigma}^m \rw T_r$.
If such $m$ exists then $T$ is called an \emph{$m$--proper obstructor}.

\end{definition}

\indent
The first example of proper obstructor is given by the following proposition.

\begin{proposition}
\label{mani-obst}
Suppose that $M$ is a $k$--dimensional closed aspherical manifold where $k>1$
and $X$ is the universal cover of $M$.    Suppose also that $X$ has a triangulation
so that $X$ is a metric simplicial complex and a group $G\!=\! \pi_1(M)$ acts on $X$
properly discontinuously, cocompactly, simplicially, and freely by isometries.
Then $X^k$ is a $(k\!-\!1)$--proper obstructor.
\end{proposition}

We prove Proposition~\ref{mani-obst} in Section~\ref{ca}.
The key ideas are the following:

\begin{enumerate}

\item Since $G$ acts on $X$ properly
discontinuously, cocompactly, simplicially, and freely by isometries,
$X$ is \emph{uniformly contractible}.  Recall that a metric space $Y$ is
\emph{uniformly contractible} if for any $r>0$, there exists
$R>r$ such that $B_r(y)$ is contractible in $B_R(y)$ for any $y \in Y$.

\item For any $R> 0$, there exists $R' > R$ so that the inclusion induced map
\[i_* \colon\thinspace \tilde{H}_j(X_{R'}) \rw \tilde{H}_j(X_R)\] is trivial for $j \! \neq \! k-1$
and $\zt \! \cong \! i_*(\tilde{H}_{k-1}(X_{R'})) \leq \tilde{H}_{k-1}(X_R)$.
(See Lemma~\ref{key-lemma}.)

\item
\label{delta}
We recall the definition of \emph{$\Delta$--complex} and use it to
complete the proof as sketched below.
\end{enumerate}

\begin{definition}
A $\Delta$--complex is a quotient space of a collection of
disjoint simplices of various dimensions, obtained by identifying some
of their faces by the canonical linear homeomorphisms that preserve
the ordering of vertices.
\end{definition}

\indent
Suppose $(\tilde{\Sigma}^m, a)$
is an essential $\zt\! - \! m \! - \!$cycle with a $\zt$--equivariant map \linebreak
$f \colon\thinspace (\tilde{\Sigma}^m, a) \rw T_r$.
Let \[|\tilde{\Sigma}^m| = \cup_{i=1}^n \Delta^m_i\]
(union of $n$--copies of $m$--simplices, use subscripts to denote different copies of
$m$--simplices) and \[ f_i \equiv f|_{\Delta^m_i}\]
Then condition (i) of Definition~\ref{ess-cycle} implies that
$\Sigma_{i=1}^n  f_i$ is an $m$--cycle of $T_r$ (over $\zt$).
That is, an essential $\zt$--$m$--cycle $(\tilde{\Sigma}^m, a)$
with a $\zt$--equivariant map $f \colon\thinspace (\tilde{\Sigma}^m, a) \rw T_r$ can
be considered as an $m$--cycle of $T_r$ (over $\zt$).
Next suppose that $g\!=\! \sum_{i=1}^n g_i$ is an $m$--chain of $T_r$
(over $\zt$) where $g_i \colon\thinspace \Delta^m \rw T_r$ are singular $m$--simplices.
Take an $m$--simplex for each $i$ and index them as $\Delta_i^m$.
Let $\Delta^{m-1}_i$ denote a codimension $1$ face of $\Delta_i^m$.
Construct a $\Delta$--complex $\Pi$ as follows:

\begin{itemize}
\item $|\Pi|=\cup_{i=1}^n \Delta^m_i$
\item For each $\ell \! \neq \! j$ we identify $\Delta^m_{\ell}$ with
$\Delta^m_j$ along $\Delta^{m-1}_{\ell}$ and $\Delta^{m-1}_j$ whenever
$g_l|_{\Delta^{m-1}_{\ell}}=g_j|_{\Delta^{m-1}_j}$.
\end{itemize}

Subdivide $\Pi$ if necessary so that $\Pi$ becomes a simplicial complex.
Consider when $g$ is an $m$--cycle and an $m$--boundary.

\indent
First, suppose $g$ is an $m$--cycle.  Then for any codimension $1$
face $\Delta^{m-1}_i$ of $\Delta^m_i$ there are an even
number of $j$'s(including $i$ itself) between $1$ and $n$ such that
$g|_{\Delta^{m-1}_i}=g|_{\Delta^{m-1}_j}$\   So $\Pi$ satisfies condition (i) of
Definition~\ref{ess-cycle} and we can consider $g$ as a map \[g \colon\thinspace \Pi \rw T_r\]
by setting $g|_{\Delta^m_i}=g_i$.

\indent
Second, suppose $g$ is an $m$--boundary.  Then there is an $(m+1)$--chain
$G \equiv \sum_{i=1}^N G_i$  where $G_i \colon\thinspace \Delta^{m+1} \rw T_r$ are singular
$(m+1)$--simplices such that $\partial G\! = \! g$.  As before one can construct a
simplicial complex $\Omega$ and consider $G$ as a map
\[ G \colon\thinspace \Omega \rw T_r \]  Let
$\partial \Omega \equiv \cup \{ m-$simplices of $\Omega$
which are the faces of an odd number of $(m\! +\! 1)$--simplices$\}$.
Note that $\partial \Omega \overset{comb.}{\cong}  \Pi$
where $\overset{comb.}{\cong}$ denotes combinatorial equivalence.
This observation will be used to construct an essential cycle in the proof
of Proposition~\ref{mani-obst}.

\section{Coarse Alexander duality}
\label{ca}

We first review the terminology of ~\cite{KK}.  Some terminology
already defined is modified in the PL category. Let $X$ be (the
geometric realization of) a locally finite simplicial complex. We
equip the $1$-skeleton $X^{(1)}$ with path metric by defining each
edge to have unit length. We call such an $X$ with the metric on
$X^{(1)}$ a metric simplicial complex. We say that $X$ has
\emph{bounded geometry} if all links have a uniformly bounded
number of simplices.  Recall that $X_r \equiv Cls(\, X^2-
N_r(d(X))\,)$, see Definition~\ref{nbd}. Also denote:
\[ \left \{ \begin{array}{l}
             B_r(c) \equiv \{x \in X |d(c,x)\leq r \} \\
             \partial B_r(c) \equiv \{x \in X |d(c,x) = r \}
                         \end{array}
   \right. \]
If $C_*(X)$ is the simplicial chain complex and $A \subset C_*(X)$ then the
\emph{support of} $A$, denoted by $Support(A)$, is the smallest subcomplex of
$K \subset X$ such that $A \subset C_*(K)$.  We say that a homomorphism
\[ h \colon\thinspace  C_*(X) \rw C_*(X)\] is \emph{coarse Lipschitz} if for each simplex
$\sigma \subset X$, $Support( h(C_*(\sigma)) )$ has uniformly bounded diameter.
We call a  coarse Lipschitz map with
 \[ D \equiv max_{\sigma} diam(\ Support( h(C_*(\sigma)) )\    ) \]
$D$-\emph{Lipschitz}.
We call a homomorphism $h$ \emph{uniformly proper}, if it is coarse Lipschitz and
there exists a proper function $\phi \colon\thinspace \rr_+ \rw \rr$ so that for each subcomplex
$K \subset X$ of diameter $\geq r$, $Support( h(C_*(\sigma)) )$ has diameter
$\geq \phi(r)$.  We say that a homomorphism $h$ has
\emph{displacement} $\leq D$ if for every simplex $\sigma \subset X$,
$Support(\ h(C_*(\sigma))\ ) \subset N_D(\sigma)$.  A metric simplicial complex is
\emph{uniformly acyclic} if for every $R_1$ there is an $R_2$ such that for each
subcomplex $K \subset X$ of diameter $\leq R_1$ the inclusion $K \rw N_{R_2}(K)$
induces zero on reduced homology groups.

\begin{definition}[(PD group)]
A group $\Gamma$ is called an $n$-dimensional \emph{Poincar\'{e} duality group}
($PD(n)$ group in short) if the following conditions are satisfied:
\begin{itemize}
\item[(i)] $\Gamma$ is of type $FP$ and $n=dim(\Gamma)$.
\item[(ii)] $H^j(\Gamma;\mathbb{Z}\Gamma) = \left \{ \begin{array}{ll}
                                          0          & j \neq n \\
                                          \mathbb{Z} & j=n
                                          \end{array}
\right.$
\end{itemize}
\end{definition}

\begin{example}
The fundamental group of a closed aspherical $k$--manifold is a $PD(k)$ group.
See~\cite{Br} for details.
\end{example}

\begin{definition}[(Coarse Poincar\'{e} duality space~\cite{KK})]
\label{pd}
A \emph{Coarse Poincar\'{e} duality space of formal dimension k},
$PD(k)$ space in short, is a bounded geometry metric simplicial complex $X$
so that $C_*(X)$ is uniformly acyclic, and there is a constant $D_0$
and chain mappings
\[C_*(X) \stackrel{\bar{P}}{\rw} C_c^{k-*}(X) \stackrel{P}{\rw} C_*(X) \]
so that
\begin{itemize}

\item[(i)] $P$ and $\bar{P}$ have displacement $\leq D_0$,
\item[(ii)] $P \circ \bar{P}$ and $\bar{P} \circ P$ are chain homotopic to the identity
by $D_0$-Lipschitz chain homotopies
$\Phi \colon\thinspace C_{*}(X) \rw C_{*+1}(X),\ \bar{\Phi} \colon\thinspace C^*_c(X) \rw C^{*-1}_c(X)$.
We call coarse Poincare duality spaces of formal dimension k a
\emph{coarse PD(k) spaces}.
\end{itemize}
\end{definition}

\begin{example}
An acyclic metric simplicial complex that admits a free, simplicial cocompact
action by a $PD(k)$ group is a coarse $PD(k)$ space.
\end{example}

\textsl{For the rest of the paper, let $X$ denote the universal cover of a
$k$-dimensional closed aspherical manifold where $k >1$}. 

Assume also that $X$ has a triangulation so that $X$ is a metric
simplicial complex with bounded geometry, and $G=\pi_1(M)$ acts on $X$
properly discontinuously, cocompactly, simplicially, and freely by
isometries. In particular, $G$ is a $PD(k)$ group and $X$ is a coarse
$PD(k)$ space. The following theorem was proved
in~\cite{KK}. Pro-Category theory is reviewed in Appendix A.

\begin{theorem}[Coarse Alexander duality~\cite{KK}]
\label{alex}
Suppose $Y$ is a coarse $PD(n)$ space, $Y'$ is a bounded geometry, uniformly acyclic
metric simplicial complex, and $f \colon\thinspace  C_*(Y') \rw C_*(Y)$ is a uniformly proper chain map.
Let $K \equiv Support(f(C_*(Y'))),\ Y_R \equiv Cls(\,Y-N_R(K)\,)$.
Then we can choose $0< r_1 < r_2 < r_3 < \ldots$ and define the inverse system
$pro\tilde{H}_j (Y_r) \equiv \{\tilde{H}_j (Y_{r_i}), i_*, \nn \}$ so that
\[pro\tilde{H}_{n-j-1}(Y_r) \cong H^j_c(Y'). \]
\end{theorem}

We rephrase the coarse Alexander duality theorem.
\begin{lemma}
\label{key-lemma}
Recall that $X$ is a metric simplicial complex with bounded geometry and a group $G$
acts on $X$ properly discontinuously, cocompactly, simplicially, and freely by isometries.
Also recall that $X_r \equiv Cls(\, X^2- N_r(d(X))\, )$.
One can choose $ 0< r_1 < r_2 < r_3 < \ldots$ and define the inverse system
$pro\tilde{H}_j (X_r)$ $\equiv \{\tilde{H}_j (X_{r_i}), i_*, \nn \}$ so that:
\[pro\tilde{H}_j (X_r) =\left \{\begin{array}{ll}
                                \mathbf{0}, & j \neq k-1 \\
                                \zt ,& j=k-1
                                 \end{array}
\right.  \]

\begin{proof}
Consider the diagonal map
\[d \colon\thinspace X \rw X^2,\ x \mapsto  (x,x) \]  and
note that $d$ is uniformly proper and $X^2$ is a $PD(2k)$ space.
Theorem~\ref{alex} implies that \[ pro\tilde{H}_{2k- \ast-1} (X_r) =
H^{\ast}_c(X). \]  Finally
observe that $H^*_c(X) \cong H_{k-*}(\rr^k) \cong H^*_c(\rr^k)$.
\end{proof}
\end{lemma}

Now we prove Proposition~\ref{mani-obst}.

\begin{figure}[ht!]\small
\begin{center}
\psfrag{sigma}{$\sigma$}
\psfrag{sigma_prime}{$\sigma '$}
\includegraphics[width=4cm]{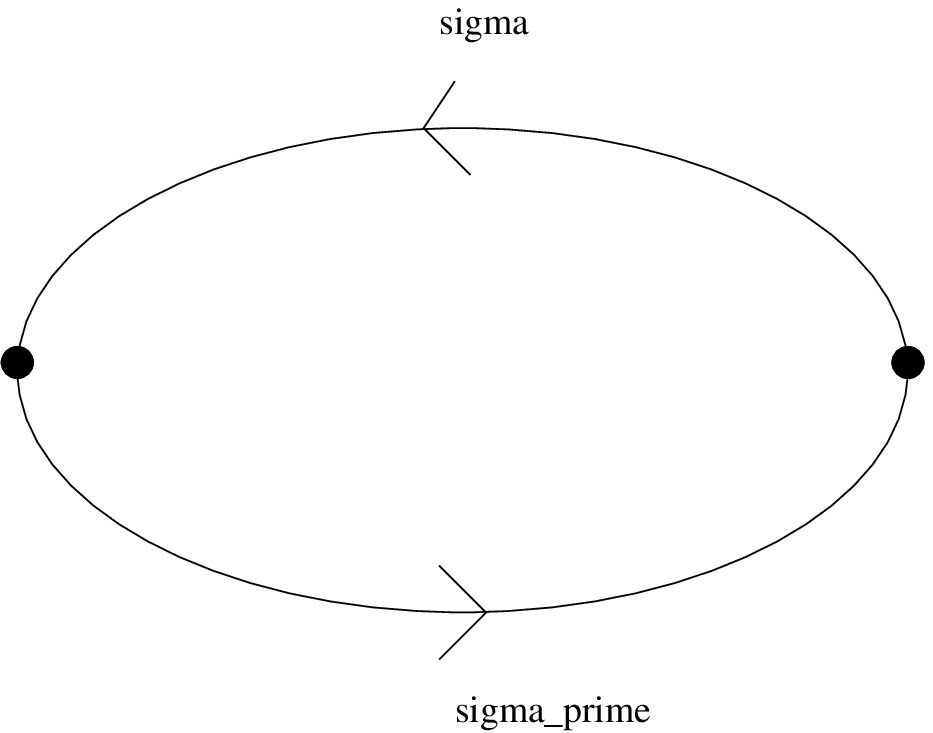}
\caption{$\sigma_1$}
\label{pic2}
\end{center}
\end{figure}

\medskip
{\bf Proof of Proposition~\ref{mani-obst}}\qua
Let $r>0$ be given.
First use Lemma~\ref{key-lemma} to choose $r= r_1 < r_2 < \ldots  < r_{k-1}<r_k$
so that \[ i_* \colon\thinspace \tilde{H}_j(X_{r_{m+1}}) \rw \tilde{H}_j(X_{r_m}) \]
is trivial for $j \neq k-1$.  In particular,
$i \colon\thinspace X_{r_k} \rw X_{r_{k-1}}$ is trivial in $\pi_0$.
Let $S^0 \equiv \{e, w\}$ and define
an involution $a_0$ on $S^0$ by $a_0(e)=w$ and $a_0(w)=e$.
Let $\theta\colon\thinspace (S^0, a_0) \rightarrow (X_{r_k}, s)$ be an equivariant map
where $s$ is the obvious involution on $X_{r_i}$.
Now let \[ \sigma \colon\thinspace I \rw X_{r_{k-1}}\] be a path so that
$\sigma(0)=\theta(e)$ and $\sigma(1)=\theta(w)$.  Define
\[ \sigma' \colon\thinspace I \rw X_{r_{k-1}}\] by $\sigma'(t)=s \sigma(t)$.  Observe that
$\sigma_1 \equiv \sigma+\sigma'$ is an $1$--cycle in $X_{r_{k-1}}$.
See Figure \ref{pic2}.

\begin{figure}[hb!]\small
\begin{center}
\psfrag{tstp}{$\tilde{\Sigma}^2_+$}
\psfrag{so}{$S^1$}
\includegraphics[width=4cm]{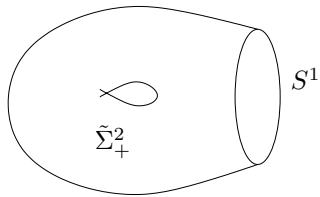}
\caption{$\tilde{\Sigma}^{2}_+$ } \label{pic3}
\end{center}
\end{figure}

Let $a_1$ be the obvious involution on $S^1$ and consider $\sigma_1$ as an
equivariant map \[ \sigma_1\colon\thinspace (S^1, a_1) \rw (X_{r_{k-1}},s) .\]
Since $i_* \colon\thinspace \tilde{H}_1(X_{r_{k-1}}) \rw \tilde{H}_1(X_{r_{k-2}})$ is trivial,
$\sigma_1$ is the boundary of a $2$--chain in $X_{r_{k-2}}$.
Call this $2$ --chain $\sigma_2^+=\sum_{i=1}^m g_i$ where
$g_i$ are singular $2$--simplices.
Following Remark~(\ref{delta}) after Proposition~\ref{mani-obst}, construct
a simplicial complex $\tilde{\Sigma}^{2}_+$ such that
\[\sigma_2^+ \colon\thinspace \tilde{\Sigma}^{2}_+ \rw  X_{r_{k-2}}\ \mbox{\ and\ }\
\partial \sigma_2^+ =\sigma_1. \]  See Figure \ref{pic3}.
Define the \emph{boundary} of $\tilde{\Sigma}^{2}_+$, $\partial
\tilde{\Sigma}^{2}_+$, to be the union of $1$--simplices, which are
the faces of an odd number of $2$--simplices.  Recall also from
Remark~(\ref{delta}) that $\partial \tilde{\Sigma}^{2}_+
\overset{comb.}{\cong} S^1$ where $\overset{comb.}{\cong}$ denotes
combinatorial equivalence.

\begin{figure}[ht!]\small
\begin{center}
\psfrag{tstp}{$\tilde{\Sigma}^2_+$}
\psfrag{x}{$x$}
\psfrag{tstm}{$\tilde{\Sigma}^2_-$}
\psfrag{aox}{$a(x)$}
\includegraphics[width=9cm]{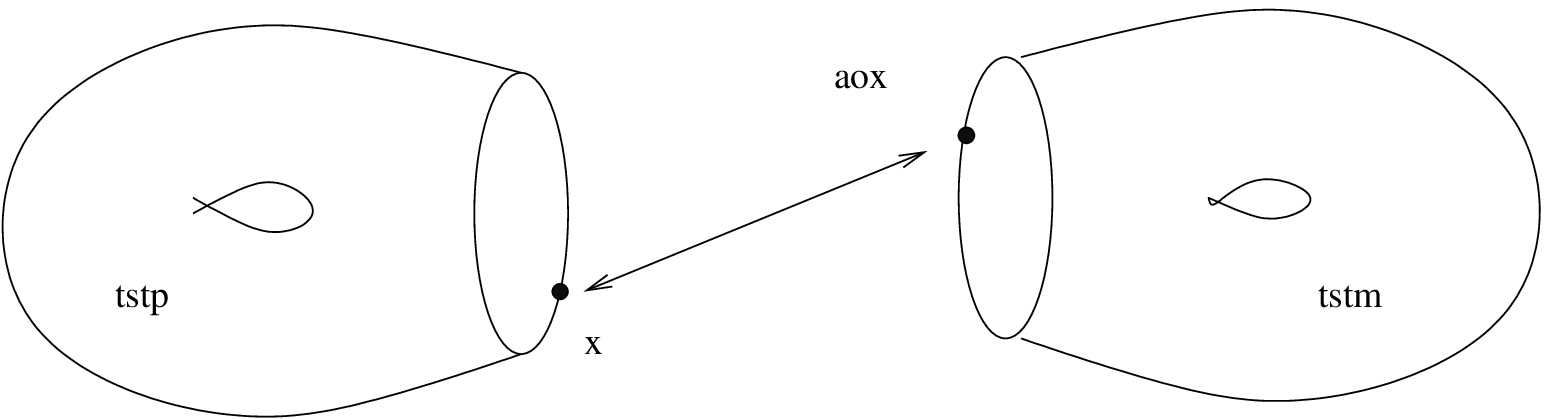}
\caption{Constructing $\tilde{\Sigma}^{2}$}
\label{pic4}
\end{center}
\end{figure}

Next, let  $\sigma_2^-=s \sigma_2^+=\sum_{i=1}^m sg_i$.  Take a copy of
$\tilde{\Sigma}^{2}_+$, denoted by $\tilde{\Sigma}^{2}_-$, such that
\[ \sigma_2^- \colon\thinspace \tilde{\Sigma}^{2}_- \rw  X_{r_{k-2}}\ \mbox{\ and\ }\
\partial \sigma_{2}^-=\sigma_1.\]
Construct $\tilde{\Sigma}^2$ by attaching $\tilde{\Sigma}^{2}_+$ and $\tilde{\Sigma}^{2}_-$
along $S^1=\partial \tilde{\Sigma}^{2}_+ =\partial \tilde{\Sigma}^{2}_- $ by identifying
$x \sim a_1(x)$.  That is, $\tilde{\Sigma}^2 \equiv \tilde{\Sigma}^{2}_+ \cup_{S^1}
\tilde{\Sigma}^{2}_-$.  See Figure \ref{pic4}.
\noindent
Define an involution $a_2$ on $\tilde{\Sigma}^2$ by setting
\[ a_2(x)=\left \{ \begin{array}{ll}
   x \in \tilde{\Sigma}^{2}_+  & \mbox{\ if\ }\  x \in \tilde{\Sigma}^{2}_- - S^1 \\
   x \in \tilde{\Sigma}^{2}_-  & \mbox{\ if\ }\  x \in \tilde{\Sigma}^{2}_+ - S^1 \\
   a_1(x) & \mbox{\ if\ }\  x \in S^1
           \end{array}
\right. \]
\noindent
Observe that $\sigma_2 \equiv \sigma_2^+ + \sigma_2^-$ is a $2$--cycle
in $X_{r_2}$ and we can consider $\sigma_2$ as an equivariant map
\[\sigma_2 \colon\thinspace (\tilde{\Sigma}^{2}, a_2) \rw  (X_{r_{k-2}}, s). \]

Continue inductively and construct a $(k\! -\! 1)$--cycle
\[ \sigma_{k-1} \colon\thinspace (\tilde{\Sigma}^{k-1}, \ a_{k-1}) \rw 
(X_{r_1}\! = \! X_{r_{k-(k-1)}}, \ s) ' \]
Simply write $a$ instead of $a_{k-1}$, and note that $X_{r_1} \subset X_r$.
So $(\tilde{\Sigma}^{k-1}, \ a)$ satisfies conditions (i)--(ii) of
Definition~\ref{ess-cycle} and we only need to show that
it satisfies condition (iii).

It was proved in~\cite{BKK} that there exists a $\zt$--equivariant
homotopy equivalence $\tilde{h} \colon\thinspace X_0 \rw S^{k-1}$. So $\tilde{h}$
induces a homotopy equivalence
\[h \colon\thinspace X_0/\sim \rw \mathbb{R}P^{k-1} .\]
Let $g \equiv h i \sigma_{k-1}  \colon\thinspace \tilde{\Sigma}^{k-1} \stackrel{\sigma_{k-1}}{\rw}
X_r \stackrel{i}{\rw} X_0 \stackrel{h}{\rw} S^{k-1}$.   Note that $g$ is equivariant.
We shall prove that $deg(g)= 1 (mod\  2)$ by constructing another map
\[f_{k-1} \colon\thinspace \tilde{\Sigma}^{k-1} \rw S^{k-1} \] with odd degree and applying
Proposition~\ref{well-prop}.

Observe that
\[S^1 \subset \tilde{\Sigma}^{2} \subset \tilde{\Sigma}^{3} \subset \ldots
\subset \tilde{\Sigma}^{k-2} \subset \tilde{\Sigma}^{k-1} \]
and for each $i$, $2 \leq i \leq k-1$:
\[\tilde{\Sigma}^{i} =\tilde{\Sigma}^{i}_+ \cup_{\tilde{\Sigma}^{i-1}}
 \tilde{\Sigma}^{i}_- \]
\noindent
Now construct a map $f_{k-1} \colon\thinspace \tilde{\Sigma}^{k-1} \rw S^{k-1}$
as follows: 
First let $f_1 \colon\thinspace S^1 \rw S^1$ be the identity and extend $f_1$ to
$f^+_2 \colon\thinspace \tilde{\Sigma}^{2}_+ \rw B^2$ by Tietze Extension theorem.
Without loss of generality assume that $(f^+_2)^{-1}(S^1) \subset S^1 \overset{comb.}{\cong}
\partial \tilde{\Sigma}^2_+$.
Then extend equivariantly to $f_2\colon\thinspace\tilde{\Sigma}^{2} \rw S^2$.
Note that $f_2^{-1}(B^2_+) \subset \tilde{\Sigma}^{2}_+$,
$f_2^{-1}(B^2_-) \subset \tilde{\Sigma}^{2}_-$, and $f_2^{-1}(S^1) \subset S^1$.

Continue inductively and construct an equivariant map
\[f_{k-1}\colon\thinspace\tilde{\Sigma}^{k-1} \rw S^{k-1}. \]
By construction, we know that \[f_j^{-1}(B^j_+) \subset \tilde{\Sigma}^{j}_+,\
f_j^{-1}(B^j_-) \subset \tilde{\Sigma}^{j}_-,\  \mbox{\ and\ }\
f_j^{-1}(S^{j-1}) \subset \tilde{\Sigma}^{j-1},\ 2 \leq j \leq k-1 .\]
Observe that $deg(f_{k-1})=deg(f_{k-2})= \ldots = deg(f_2)=deg(f_1)$.
(Recall that $deg(f_m) \equiv$\ the number of $m$--simplices of
$\tilde{\Sigma}^m$ mapped into a simplex $s$ of $S^m$ by $f$.)
But $deg(f_1)=id_{S^1}=1(mod\ 2)$
so  $f_{k-1}\colon\thinspace \tilde{\Sigma}^{k-1} \rw S^{k-1}$
has nonzero degree.  Now Proposition~\ref{well-prop} implies that
$deg(g)=1(mod\ 2)$.   Therefore $(\tilde{\Sigma}^{k-1},a)$ with
equivariant map \[\sigma_{k-1}\colon\thinspace \tilde{\Sigma}^{k-1} \rw X_r .\]
satisfies conditions (i),(ii), and (iii) of Definition~\ref{ess-cycle}.
Now the proof of Proposition~\ref{mani-obst} is complete.

\section{New proper obstructors out of old}
\label{pod1}

In this Section, we construct a $k$--proper
obstructor from a $(k\! -\! 1)$--proper obstructor $X$.

\begin{definition}
\label{uray}
Let $(Y, d_Y)$ be a proper metric space and $(\alpha, d_{\alpha})$ be a
metric space isometric to $[0, \infty)$.
Let $\phi \colon\thinspace[0, \infty) \rightarrow \alpha$ be an
isometry and denote $\phi(t)$ by $\alpha_t$.  Define a metric space
$(Y \vee \alpha, d)$, called \emph{$Y$ union a ray}, as follows:

\begin{itemize}
\item[(i)] As a set $Y \vee \alpha$ is the wedge sum.  That is,
$Y \vee \alpha=Y \cup \alpha$ with $Y \cap \alpha = \{\alpha_0\}$

\item[(ii)] The metric $d$ of $Y \vee \alpha$ is defined by
\[ \left \{ \begin{array}{ll}
           d(v,w)=d_Y(v,w), & \mbox{\ if\ }\  v,w \in Y \\
           d(v,w)=d_{\alpha}(v,w), & \mbox{\ if\ }\  v,w \in \alpha\\
d(v,w)=d_Y(v, \alpha_0) + d_{\alpha}(\alpha_0, w),& \mbox{\ if\ }\  v \in Y,\  w \in \alpha
            \end{array}
\right. \]
\end{itemize}

\end{definition}

\begin{proposition}
\label{union-ray}
Let $X$ be a $k$-dimensional contractible manifold without
boundary and $k >1$.  Suppose also that $X$ has a triangulation so that $X$ is a metric
simplicial complex and a group $G$ acts on $X$ properly
discontinuously, cocompactly, simplicially, and freely by isometries.
In particular, $X$ is a $(k-1)$--proper obstructor.
Then $X \vee \alpha$ is a $k$--proper obstructor.

\begin{proof}
Recall that by Lemma~\ref{key-lemma}, we can choose
$0< r_1<r_2< r_3 \dots $ and define $pro\tilde{H}_{k-1}(X_r) \equiv
\{\tilde{H}_{k-1}(X_{r_i}), i_*, \nn \}$ so that $pro\tilde{H}_{k-1}(X_r)\! =\! \zt$.
This means that for any $r>0$ we can choose $R >r$ so that
\[ r' \geq R \Rightarrow \zt = i_*(H_{k-1}(X_{r'})) \leq H_{k-1}(X_r). \]  Now let $r>0$
be given and choose $R>r$ as above.  Let $(\tilde{\Sigma}^{k-1},a)$ be an essential
$\zt$--$(k-1)$--cycle with a $\zt$--equivariant map
\[f\colon\thinspace \tilde{\Sigma}^{k-1} \rw X_R. \]
Next consider composition $i \circ f \colon\thinspace \tilde{\Sigma}^{k-1}
\stackrel{f}{\rw} X_R \stackrel{i}{\rw} X_r$.   If $i \circ f \in Z_{k-1}(X_r)$ is
the boundary of a $k$--chain then we can construct an essential $\zt$--$k$--cycle
with $\zt$--equivariant map into $X_r$ using the method used in the
proof of Proposition~\ref{mani-obst}.   But this implies $X$ is a $k$--proper obstructor.
(Recall that $X^k$ is a $(k\! -\! 1)$--proper obstructor.)
So we can assume $i \circ f \in Z_{k-1}(X_r)\! -\! B_{k-1}(X_r)$.  That is,
$0 \neq [i \circ f] =i_* [f] \in H_{k-1}(X_r)$.  Let $p_i\colon\thinspace X_r \rw X$ denote the projection
to the $i$-th factor, $i=1,2$.

We need the following lemma.
\begin{lemma}
\label{intlem}
Define $j\colon\thinspace X-B_R \rw X_R$, $x \mapsto (\alpha_0, x)$.  Then the composition
\[ i_* \circ j_*\colon\thinspace H_{k-1}(X-B_R(\alpha_0)) \stackrel{j_*}{\rw} H_{k-1}(X_R)
\stackrel{i_*}{\rw} H_{k-1}(X_r) \]  is nontrivial.

\begin{proof}[The proof of Lemma~\ref{intlem}]
Consider a map $\lambda\colon\thinspace H_{k-1}(X_0)\rw \zt$ given by
\[ [f] \mapsto Lk(f, \Delta)(mod\ 2) \]
where $Lk(f, \Delta)$ denote the linking number
of $f$ with the diagonal $\Delta$.\footnote{We can compute $Lk(f, \Delta)$
by letting $f$ bound a chain $\tilde{f}$ transverse to $\Delta$ and setting
$Lk(f, \Delta) = Card (\tilde{f}^{-1}(\Delta))$.}
Now consider the composition:
\[\zeta \colon\thinspace H_{k-1}(X-B_R(\alpha_0)) \stackrel{j_*}{\rw} H_{k-1}(X_R) \stackrel{i_*}{\rw}
H_{k-1}(X_0) \stackrel{\lambda}{\rw} \zt \]  We shall show that $\zeta$ is
nontrivial.  Choose $[f_1] \in H_{k-1}(X-B_R)$ so that
$Lk(f_1, \alpha_0) \neq 0$ where $[\alpha_0] \in H_0(X)$.
Then $Lk(i_*j_*([f_1]), \Delta) \neq 0$.(We can choose the same chain
transverse to $\Delta$.)  Hence $\zeta$ is nontrivial.
In particular, $i_* \circ j_*$ and $j_*$ are nontrivial.
\end{proof}
\end{lemma}

Since $j_* \colon\thinspace H_{k-1}(X -B_R) \rw H_{k-1}(X_R)$ is nontrivial, we can
choose $h \in Z_{k-1}(X-B_R)\! -\! B_{k-1}(X-B_R)$ with
$g \equiv j \circ h \in Z_{k-1}(X_R)\! -\! B_{k-1}(X_R)$.  That is, $0 \neq [g] \in H_{k-1}(X_R)$.
We can consider $g$ as a map $g \colon\thinspace \Pi \rw X_R$ where
$\Pi$ is a $(k\! -\! 1)$--dimensional simplicial complex satisfying
condition (i) of Definition~\ref{ess-cycle} such that

\begin{itemize}
\item $0 \neq i_*[g] \in H_{k-1}(X_r)$
\item  $i \circ g\colon\thinspace \Pi \stackrel{g}{\rw} X_R \stackrel{i}{\rw} X_r$ with
$p_1(\  i \circ g(\Pi)\ )=\{\alpha_0\}=X \cap \alpha(0)$.
\end{itemize}

Next define $g'=sg$, that is,
\[g'\colon\thinspace \Pi \stackrel{g}{\rw} X_R \stackrel{s}{\rw} X_R. \]
Note that $i \circ g'$ is a
cycle in $X_R$ and $p_2 (\ i \circ g'(\Pi)\ )=\alpha_0$.  Also
$[f], [g] \in H_{k-1}(X_R)$ and $i_*[f],  i_*[g] \in H_{k-1}(X_r)$ are nonzero.
Observe that $i \circ f$ and $i \circ g$ must be homologous in $X_r$
since $\zt = i_*(H_{k-1}(X_R)) \leq H_{k-1}(X_r)$.
We simply write $f$, $g$,  and $g'$ instead of $i \circ f$, $i \circ g$, and $i \circ g'$.
There exists a $k$--chain $G \in C_k(X_r)$ such that
\[ \partial G =f +g .\]  Again consider $G$ as a map $G \colon\thinspace \Omega \rw X_r$
where $\Omega$ is a simplicial complex so that
\[ \partial \Omega= \tilde{\Sigma}^{k-1} \sqcup \Pi  .\]
See Figure \ref{pic5}.
\begin{figure}[ht!]\small
\begin{center}
\psfrag{pi}{$\Pi$}
\psfrag{om}{$\Omega$}
\psfrag{tskm}{$\tilde{\Sigma}^{k-1}$}
\includegraphics[width=4cm]{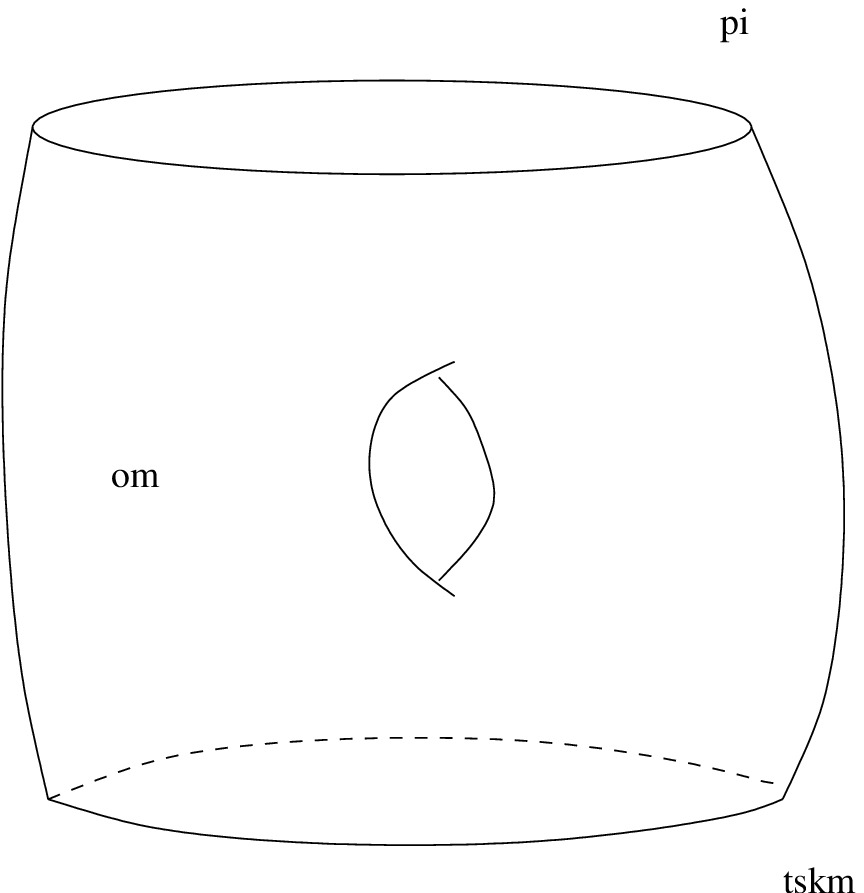}
\caption{$\Omega$}
\label{pic5}
\end{center}
\end{figure}

Next define $G'=sG$, that is,
\[G' \colon\thinspace \Omega \stackrel{G}{\rw} X_r \stackrel{s}{\rw} X_r .\]
Note that \[ \partial G' =f +g'. \]
Now take two copies of $\Omega$ and index them
as $\Omega_1$ and $\Omega_2$.  Similarly
$\Pi_1 \subset \partial \Omega_1$ and  $\Pi_2 \subset \partial \Omega_2$.
Hence \[\partial \Omega_i = \tilde{\Sigma}^{k-1} \cup \Pi_i\ , i=1,2. \]
Denote $id(x)\! =\! x'$ for $x \in \Omega_1\! -\! \tilde{\Sigma}^{k}$
where $id\colon\thinspace \Omega_1 \rw \Omega_2$.
Construct a $k$--dimensional simplicial complex
$\tilde{\Omega}$ by attaching $\Omega_1$ and $\Omega_2$
along $\tilde{\Sigma}^{k-1}$ by
$a\colon\thinspace \tilde{\Sigma}^{k-1} \rw \tilde{\Sigma}^{k-1}$.
That is,
\[\tilde{\Omega}= (\Omega_1\cup \Omega_2) / x \sim ax, \ x \in
\tilde{\Sigma}^{k-1}  .\]  See Figure \ref{pic6}.
\begin{figure}[hb!]\small
\begin{center}
\psfrag{oo}{$\Omega_1$}
\psfrag{x}{$x$}
\psfrag{ot}{$\Omega_2$}
\psfrag{aox}{$a(x)$}
 \includegraphics[height=3cm]{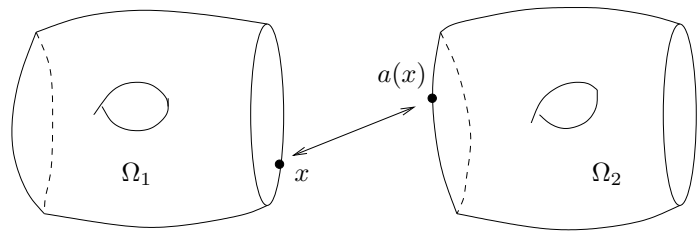}
\caption{Constructing $\tilde{\Omega}$}
\label{pic6}
\end{center}
\end{figure}

We can define an involution $\bar{a}$ on $\tilde{\Omega}$ by
\[ \left\{\begin{array}{ll}
         \bar{a}(x)=a(x), &  x \in \tilde{\Sigma}^{k-1},\\
         \bar{a}(x)=x', & x \in \Omega_1-\tilde{\Sigma}^{k-1},\\
         \bar{a}(x')=x,&  x' \in \Omega_2-\tilde{\Sigma}^{k-1}.
         \end{array}
\right. \]
Also we can define a $\zt$--equivariant map
$\Phi \colon\thinspace \tilde{\Omega} \rw X_r$ by:
\[ \left\{\begin{array}{l}
         \Phi|_{\Omega_1}=G\\
         \Phi|_{\Omega_2}=G'
          \end{array}
\right. \]

We define  \[  \tilde{\Sigma}^{k}=
(\ \Pi_1 \! \times \! [0,1] / (\Pi_1, 1) \sim *\ )
\ \cup_{\Pi_1}\  \tilde{\Omega}
\ \cup_{\Pi_2} \  (\Pi_2 \! \times \! [0,-1]
/(\Pi_2, -1) \sim *\ )    .\]   See Figure \ref{pic7}.
\begin{figure}[ht!]\small
\begin{center}
\psfrag{tsk}{$\tilde{\Sigma}^k$}
\psfrag{tskm}{$\tilde{\Sigma}^{k-1}$}
\psfrag{oo}{$\Omega_1$}
\psfrag{ot}{$\Omega_2$}
\psfrag{to}{$\tilde{\Omega}$}
\includegraphics[width=10cm]{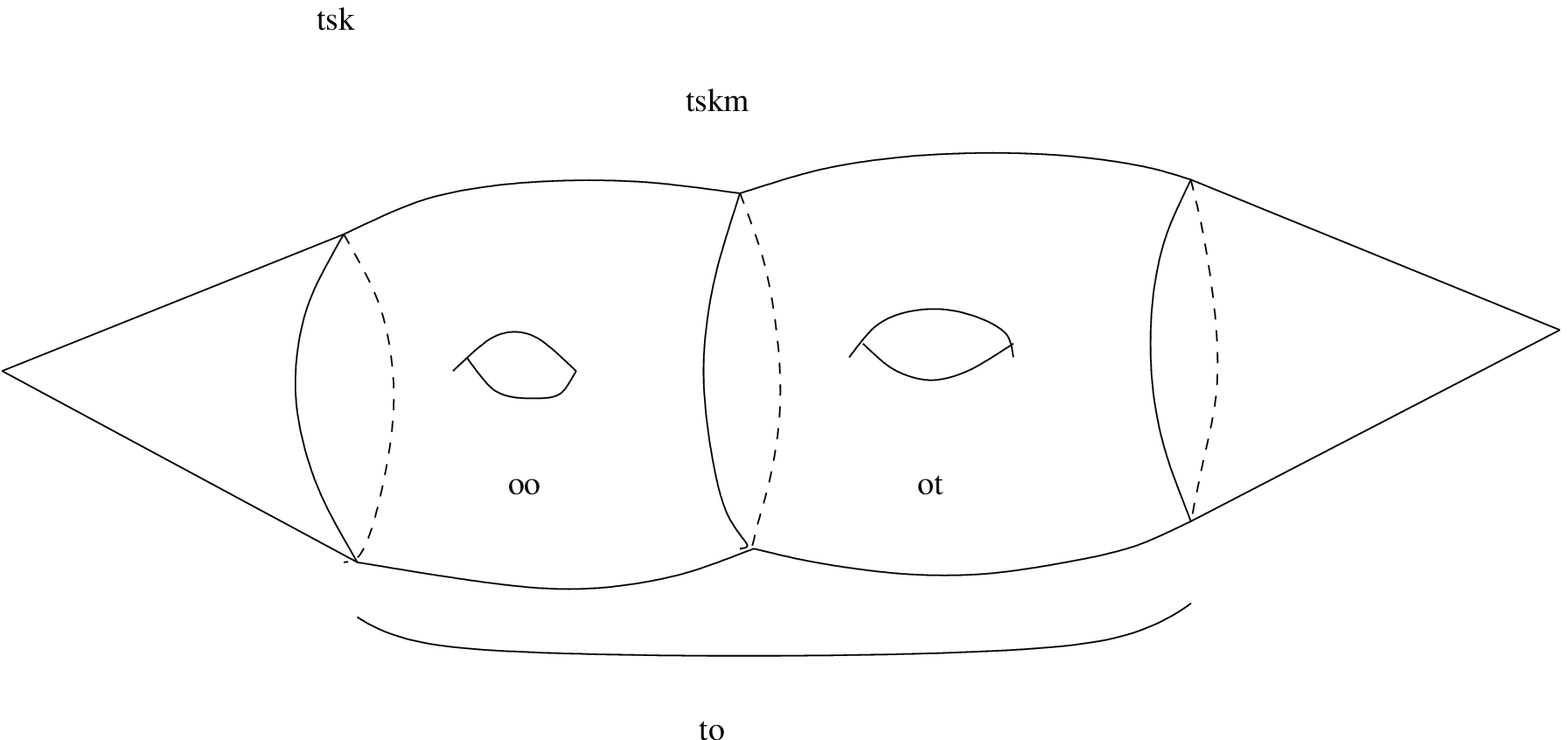}
\caption{Constructing $\tilde{\Sigma}^{k}$}
\label{pic7}
\end{center}
\end{figure}
Now extend $\bar{a}$ over $\tilde{\Sigma}^{k}$, and denote
$\tilde{\Sigma}^{k}/ x \sim a(x)$ by $\Sigma^{k}$.

Suppose that $\Sigma^{k}$ classifies into $\mathbb{R}P^m$ where $m < k$.
Let \[ h\colon\thinspace \Sigma^{k} \rw \mathbb{R}P^m \] be the classifying map
and \[\tilde{h}\colon\thinspace \tilde{\Sigma}^{k} \rw S^m \] be the equivariant map
covering $h$.  Observe that
\[deg\  \tilde{h}|_{\tilde{\Sigma}^{k-1}} =deg\  \tilde{h}=0\ (mod\ 2)\]
This is a contradiction since there already exists a
$\zt$--equivariant map
\[ \varphi\colon\thinspace \tilde{\Sigma}^{k-1} \rw S^{k-1} \]
of odd degree.  Hence $(\tilde{\Sigma}^{k}, \bar{a})$
is an essential $\zt$--$k$--cycle.

Finally, we need to define a $\zt$--equivariant map:
\[F\colon\thinspace \tilde{\Sigma}^{k} \rw (X \vee \alpha)_r \]
Recall that $p_1g(\Pi)=\alpha_0$ and let  \[c \colon\thinspace p_2g(\Pi_1)\! \times \! I \rw X \]
be a contraction to $\alpha_0$.   Similarly $p_2g'(\Pi)=\alpha_0$ and let
\[c' \colon\thinspace p_1g'(\Pi_1)\! \times \! I \rw X \] be a contraction to $\alpha_0$. \vspace{1ex}
Define a $\zt$--equivariant map
\[F\colon\thinspace \tilde{\Sigma}^{k} \rw (X \vee \alpha)_r \]
as follows: Recall that $\phi(t)=\alpha_t$ in Definition~\ref{uray}, so
$d(\alpha_0, \alpha_s)=s$ and $d(\alpha_s , x) \geq s$ for any $x \in X$.
\[ \left \{ \begin{array}{ll}
F|_{\tilde{\Omega}}=\Phi \\
F(x, t)= (\ \alpha_{2rt},\  p_2 g(x) \ ), &  x \in \Pi_1,
\ t \in [0,\frac{1}{2}] \\
F(x, t)= (\ \alpha_r,\  c_{(2t-1)}(p_2g(x)) \ ), & x \in \Pi_1,
\ t \in [\frac{1}{2},1] \\
F(x, t)= (\ p_1 g'(x), \ \alpha_{-2rt}\ ), &  x \in \Pi_2,
\ t \in [0,-\frac{1}{2}] \\
F(x, t)= (\ c_{(-2t-1)}(p_1g'(x)),\  \alpha_r\ ), &  x \in \Pi_2,
\ t \in [-\frac{1}{2},-1]
            \end{array}
\right. \]
The proof of Proposition \ref{union-ray} is now complete.
\end{proof}
\end{proposition}

If $Y$ and $Z$ are metric spaces we use the \emph{sup metric} on
$Y\!  \times \! Z$ where
\[ d_{sup}((y_1,z_1), (y_2,z_2))
\equiv max\{d_Y(y_1,y_2), d_Z(z_1,z_2) \} .\]

\begin{proposition}
\label{join-lemma}
Suppose $X_1, X_2$ are $m_1, m_2$--proper obstructors, respectively.  Then
$X_1\! \times \! X_2$ is an $(m_1\! +\! m_2\! +\! 1)$--proper obstructor.

\begin{proof}
Let $r>0$ be given and let
\[ \left\{\begin{array}{l}
          f_1\colon\thinspace \tilde{\Sigma}_{1}^{m_1} \rw (X_1)_r \\
          f_2\colon\thinspace \tilde{\Sigma}_{2}^{m_2} \rw (X_2)_r
          \end{array}
\right. \]
be $\zt$--equivariant maps for essential $\zt$--cycles.
Note that \[(X_1 \! \times \! X_2)_r =((X_1)_r\! \times \! {(X_2)}^2) \cup_{(X_1)_r \! \times \!
(X_2)_r} ({(X_1)}^2 \! \times \! (X_2)_r).\]  Let $a_1$ be the involution on $(X_1)_r$
and $a_2$ be the involution on $(X_2)_r$.  Recall that the join
$\tilde{\Sigma}_{1}^{m_1} * \tilde{\Sigma}_{2}^{m_2}$ is obtained from
$\tilde{\Sigma}_{1}^{m_1} \! \times \! \tilde{\Sigma}_{2}^{m_2} \! \times \! [-1,1]$ by
identifying $\tilde{\Sigma}_{1}^{m_1} \! \times \! \{y\} \! \times \! \{1\}$ to a point
for every $y \in  \tilde{\Sigma}_{2}^{m_2}$ and identifying
$\{x\} \! \times \! \tilde{\Sigma}_{2}^{m_2} \! \times \! \{-1\}$ to a point for every
$x \in  \tilde{\Sigma}_{1}^{m_1}$.  Define an involution $a$ on
$\tilde{\Sigma}_{1}^{m_1} * \tilde{\Sigma}_{2}^{m_2}$ by
\[a(v,w,t)=(a_1(v), a_2(w),t).\]
Let
\[ \left\{\begin{array}{l}
           g_1\colon\thinspace \tilde{\Sigma}_1^{m_1} \rw S^{m_1}\\
          g_2\colon\thinspace \tilde{\Sigma}_2^{m_2} \rw S^{m_2}
          \end{array}
\right. \]
be equivariant maps of odd degree.  Then:
\begin{gather*}
g_1 * g_2 \colon\thinspace \tilde{\Sigma}_{1}^{m_1} * \tilde{\Sigma}_{2}^{m_2}
\rw S^{m_1}*S^{m_2} = S^{m_1+m_2+1}\\
(v,w,t) \mapsto (g_1(v), g_2(w), t)\end{gather*}
is also an equivariant map of an odd degree. Hence
$(\tilde{\Sigma}_{1}^{m_1} * \tilde{\Sigma}_{2}^{m_2},a)$
is an essential $\zt$--$(m_1\! +\! m_2\! +\! 1)$--cycle.

Now let
\[c\colon\thinspace f_1(\tilde{\Sigma}_{1}^{m_1}) \! \times \! [-1,1] \rw X_1^2 \]
be a $\zt$--equivariant contraction to a point  such that
$c_t=id$ for $t \in [-1, 0]$.  Similarly let
\[d\colon\thinspace f_2(\tilde{\Sigma}_{2}^{m_2}) \! \times \! [-1,1] \rw X_2^2 \]
be a $\zt$--equivariant contraction to a point such that $d_t=id$
for $t \in [0, 1]$.

Finally define
\[ f\colon\thinspace \tilde{\Sigma}_{1}^{m_1} * \tilde{\Sigma}_{2}^{m_2} \rw
(X_1 \! \times \! X_2)_r\ \mbox{by}\  f(v,w,t)=(\, c_t(f_1(v)), d_t(f_2(w))\, )
. \]
We note that $f$ is $\zt$--equivariant.
\end{proof}
\end{proposition}

\section{Proper obstructor dimension}
\label{pod2}

We review one more notion from~\cite{BKK}.
\begin{definition}
The \emph{uniformly proper dimension}, $updim(G)$, of a discrete
group $G$ is the smallest integer $m$ such that there is a
contractible $m$--manifold $W$ equipped with a proper metric $d_W$,
and there is a $g\colon\thinspace \Gamma \rw W$ with the following properties:
\begin{itemize}
\item $g$ is Lipschitz and uniformly proper.
\item There is a function $\rho\colon\thinspace (0, \infty)\rw (0,\infty)$ such that any
ball of radius $r$ centered at a point of the image of $h$ is contractible in the ball
of radius $\rho(r)$ centered at the same point.
\end{itemize}
If no such $n$ exists, we define $updim(G)=\infty$.
\end{definition}
\noindent
It was proved in~\cite{BKK} that
\[ updim(G) \leq actdim(G). \]

Now we generalize the obstructor dimension of a group.

\begin{definition}
\label{pobdim}
The \emph{proper obstructor dimension} of $G$, $pobdim(G)$,
is defined to be $0$ for finite groups,\ $1$ for $2$--ended
groups, and otherwise $m\! +\! 1$ where $m$ is the largest integer
such that for some $m$--proper obstructor $Y$, there exists a uniformly proper
map \[ \phi\colon\thinspace Y \rw T_G\] where $T_G$ is a proper metric space with a quasi-isometry
$q \colon\thinspace T_G \rw G$.
\end{definition}

\begin{lemma}
Let $Y$ be an $m$--proper obstructor.  If there is a uniformly proper
map $\beta\colon\thinspace Y \rw W^d$ where $W$ is a contractible $d$--manifold then
$d >m$.

\begin{proof}
Assume $d\!  \leq \! m$ $(d\! -\! 1 \leq m\! -\! 1)$.
Observe that if $\beta$ is uniformly proper then $\beta$ induces
an equivariant map $\bar{\beta}\colon\thinspace Y_r \rw W_0$ for some large $r>0$.
Now let $f\colon\thinspace \tilde{\Sigma}^{m} \rw Y_r$ be an essential
$\zt\! -\! m \! - \!$cycle where $f$ is equivariant.  Let $h\colon\thinspace W_0 \rw S^{d-1}$
be an equivariant homotopy equivalence.  We have an equivariant map
\[g=ih\bar{\beta}f\colon\thinspace \tilde{\Sigma}^{m} \stackrel{f}{\rw}  Y_r
\stackrel{\bar{\beta}}{\rw} W_0 \stackrel{h}{\rw}
S^{d-1} \stackrel{i}{\rw} S^{m-1} \stackrel{i}{\rw} S^m \]
where $i\colon\thinspace S^{d-1} \rw S^{m-1} \rw S^m$ is the inclusion. 
Note that $g$ is equivariant but $deg(g)\! =\! 0 (mod \ 2)$.  This is a
contradiction by Proposition~\ref{well-prop}.
\end{proof}

\end{lemma}

Suppose that $G$ is finite so that $pobdim(G)=0$ by definition.

Clearly, $actdim(G)= 0$ if $G$ is finite.  Hence $pobdim(G)\! =\! 
actdim(G)\! =\! 0$ in this case.   Next suppose that $G$ has two ends so that
$pobdim(G)\! =\! 1$.  Note that there exists $\zz \cong H \leq G$ with
$|G \! :\!  H|\!  <\!  \infty$.  And this implies that
\[ actdim(G) = actdim(H)=actdim (\zz)=1 .\]
Therefore, $pobdim(G)\! =\! actdim(G)\! =\! 1$ when $G$ has two ends.
Now we prove the main theorem for the general case.

\begin{maint}
$pobdim(G) \leq updim(G)\leq actdim(G)$

\begin{proof}
We only need to show the first inequality.   Let $pobdim(G)\! =\! m\! +\! 1$ for some $m\!  >\! 0$.
That is, there exists an $m$--proper obstructor $Y$, a proper metric space $T_G$,
a uniformly proper map $\psi\colon\thinspace Y \rw T_G$,  
and a quasi-isometry $q\colon\thinspace T_G \rw G$.
Let  $updim(G) \! \equiv \! d$ such that there exists a
uniformly proper  map $\phi \colon\thinspace G \rw W^d$ where $W$ is a contractible $d$--manifold.
But the composition
\[ \phi \! \circ \! q \! \circ \! \psi\colon\thinspace Y \rw T_G \rw G \rw W^d \]
is uniformly proper.  Therefore
\[ m+1=pobdim(G) \leq updim(G) \]
by the previous lemma.
\end{proof}

\end{maint}

Before we consider some applications, we make the following observation about
compact aspherical manifolds with incompressible boundary.

\begin{lemma}
Assume that $W$ is a compact aspherical $m$--manifold with all boundary components
incompressible.  Let $\pi \colon\thinspace \tilde{W} \rw W$ denote the universal cover of $W$.
Suppose that there is a component of $\partial W$,   call it $N$, so that
$|\pi_1(W) \colon\! \pi_1(N)| >\! 2$.  
Then  $|\pi_1(W) \colon \! \pi_1(N)|$ is infinite.

\proof
Observe that $N$ is aspherical also.   First, we show that if
\[ 1< |\pi_1(W)\colon \! \pi_1(N)| < \infty \]
then $\tilde{M} \equiv \tilde{W}/ \pi_1(N)$ has two
boundary components and $W$ has one boundary component.  We claim that $\tilde{M}$
has a boundary component homeomorphic to $N$ which is still denoted by $N$.
To see this consider the long exact sequence:
\[ \cdots \rw H_1(\partial \tilde{M}) \stackrel{i_*}{\rw} H_1(\tilde{M}) \rw
H_1(\tilde{M}, \partial \tilde{M}) \rw \tilde{H}_0(\partial \tilde{M}) \rw \tilde{H}_0(\tilde{M}) =0 \]
Since $\pi_1(N)=\pi_1(\tilde{M})$,  $i_*\colon\thinspace H_1(\partial \tilde{M}) \rw H_1(\tilde{M})$  is surjective.
So we have:
\[ 0  \rw H_1(\tilde{M}, \partial \tilde{M}) \rw \tilde{H}_0(\partial \tilde{M}) \rw 0\]
Since $|\pi_1(W)\colon\! \pi_1(N)| $ is finite $\tilde{M}$ is compact.
Now $H_1(\tilde{M}, \partial \tilde{M}) \cong H^{m-1} (\tilde{M})$ by duality.
But $H^{m-1} (\tilde{M}) \cong H^{m-1} (N)$ and $H^{m-1} (N) \cong \zt$ since $N$ is a closed
$(m\! -\! 1)$--manifold.    That is, $\tilde{H}_0(\partial \tilde{M}) \cong \zt$ so $\tilde{M}$ has two
boundary components.  Next
let $N$ and $N'$ denote two boundary components of $\partial \tilde{M}$  both of which are
mapped to $N \subset W$ by $p\colon\thinspace \tilde{M} \rw W$.   Hence $\partial W$ has one component.

Now assume that $m \equiv |\pi_1(W)\colon\! \pi_1(N)| >\! 2$.  Suppose $m$ is finite.
Note that  $p|_{N} \colon\thinspace N (\subset \tilde{M})  \rw N (\subset W)$ has index $1$,
and $p|_{N'} \colon\thinspace N' (\subset \tilde{M})  \rw N (\subset W)$ has index $m-1$.
This means that $|\pi_1(\tilde{M}) \colon\! \pi_1(N')|=m\! -\! 1$  since
$\pi_1(\tilde{M} ) =\pi_1(N)$. There are two alternative arguments:
\begin{itemize}
\item If $m \! >\! 2$ then $\tilde{M}$ is an aspherical manifold with two boundary components $N$
and $N'$ with $|\pi_1(\tilde{M}) \colon\! \pi_1(N')|=m\! -\! 1>\! 1$.   Consider $\tilde{W}/\pi_1(N')$.
The same argument applied to $\tilde{W}/\pi_1(N')$ shows that $\tilde{M}$ has one boundary
component, which is a contradiction.  Therefore $|\pi_1(W)\colon\! \pi_1(N)|$ is infinite.
\item Suppose $m\! >\! 2$.  Choose a point $x \in\! N \! \subset \! \partial W$ and let $\tilde{x} \in\! N
\! \subset\!  \partial\!  \tilde{M}$ so that $p(\tilde{x})\! =\! x$.  Next choose two loops $\alpha$ and
$\beta$ in $W$ based at $x$ so that $\{ \pi_1(N), [\alpha] \pi_1(N),   [\beta] \pi_1(N) \}$ are
distinct cosets.
(We are assuming $|\pi_1(W)\colon \! \pi_1(N)| >\! 2$.)  Let $\tilde{\alpha}$ and $\tilde{\beta}$
be the liftings of $\alpha$ and $\beta$ respectively so that
$\tilde{\alpha}(0)\! =\! \tilde{x} \! =\! \tilde{\beta}(0)$.
Note that $\tilde{y}_1 \! \equiv \! \tilde{\alpha}(1), \tilde{y}_2 \equiv \tilde{\beta}(1) \in N'$ and
$\tilde{y}_1  \neq   \tilde{y}_2$ since $[\alpha] \pi_1(N) \neq [\beta] \pi_1(N)$.
Now consider a path $\tilde{\gamma}$ in $N'$ from $\tilde{y}_1$ to $\tilde{y}_2$.
Observe that $p\tilde{\gamma} \equiv \gamma$ is a loop based at $x$, and
$[\gamma] \in p_*(\pi_1(N')) \leq \pi_1(N)$.  But $[\alpha] [\gamma] [\beta]^{-1} =\! 1$ and
this implies that $[\alpha]^{-1}[\beta] \in \! \pi_1(N)$ contary to
$[\alpha] \pi_1(N) \neq [\beta] \pi_1(N)$.\endproof
\end{itemize}
\end{lemma}

\begin{corollary}[Application]
Suppose that $W$ is a compact aspherical $m$--manifold with incompressible boundary.
Also assume that there is a component of $\partial W$, call it $N$, so that
$|\pi_1(W)\colon\! \pi_1(N)| >\! 2$.\newline
Then $actdim(\pi_1(W))\! =\! m $.

\begin{proof}
Let $p\colon\thinspace \tilde{W} \rw W$ be the universal cover of $W$.
It is obvious that $actdim(\pi_1(W)) \leq m$ as $\pi_1(W)$ acts cocompactly and properly
discontinuously on $\tilde{W}$.   Denote $G \equiv \pi_1(W)$ and $H \equiv \pi_1(N)$.
Let $\tilde{N}$ be a component
of $p^{-1}(N)$.     Therefore $\tilde{N}$ is the contractible universal cover of $N^{(m-1)}$.
Note that $\tilde{N}$ is an $(m\! -\! 2)$--proper obstructor by Proposition~\ref{mani-obst}.
Now $\tilde{W} / H$ has a boundary component homeomorphic to $N$.  Call this
component $N$ also.    $|G \colon\! H|$ is infinite by the previous lemma, and this implies that
$\tilde{W} / H$ is not compact.   In particular,
there exists a map $\alpha' \colon\thinspace [0,\infty) \rw \tilde{W} / H$ with the
following property: For each $D >0$ there exists $T \in [0,\infty)$ such that  for any
$x \in N$,  $d(\alpha' (t), x) >D$ for $t > T$,  and $\alpha'(0) \in N$.
Let $\tilde{\alpha} \colon\thinspace [0,\infty) \rw \tilde{W}$ be a lifting of $\alpha'$ such that
$\tilde{\alpha}(0) \in \tilde{N}$.   Now we define a uniformly proper map:
\[\phi\colon\thinspace \tilde{N} \vee \alpha \rw \tilde{W}  \] 
$$\left\{\begin{array}{l}
          \phi|_{\tilde{N}}=inclusion \\
          \phi(\alpha_t)=\tilde{\alpha}(t)
          \end{array}
\right.\leqno{\mbox{}}$$

Observe that $\phi$ is a uniformly proper map.   Since $\tilde{N} \vee \alpha$ is an
$(m\! -\! 1)$--proper obstructor and $\tilde{W}$ is quasi-isometric to $G$, $pobdim(G)\! \geq \! m$.
But   \[pobdim(G) \leq updim(G) \leq actdim(G) \leq m.\]
The last inequality follows from the fact that $G$ acts on
$\tilde{W}$ properly discontinuously.  Therefore $pobdim(G) \! = \! m$.
\end{proof}
\end{corollary}

The following corollary answers Question $2$ found in \cite{BKK}.
\begin{corollary}[Application]
Suppose that $W_i$ is a compact aspherical
$m_i$--manifold with incompressible boundary for $i\! =\! 1, \ldots, d$.  Also assume
that for each $i$, $1\! \leq i \! \leq \! d$, there is a component of $\partial W_i$,
call it $N_i$, so that $|\pi_1(W_i)\colon\! \pi_1(N_i)| >\! 2$.
Let $G \equiv \pi_1(W_1) \! \times \ldots \times \! \pi_1(W_d)$.   Then:
\[actdim(G)=m_1+ \ldots + m_d \]

\proof
It is easy to see that
\[ actdim(G) \leq m_1+ \ldots + m_d \]
as $G$ acts cocompactly and properly
discontinuously on $\tilde{W}_1 \!  \times \cdots \times \! \tilde{W}_d$.
Denote $\pi_1(W_i) \equiv G_i$ and $\pi_1(N_i) \equiv H_i$.
Let \[p\colon\thinspace \tilde{W}_i \rw W_i \]
be the contractible universal cover and let $\tilde{N}_i$ be a component
of $p^{-1}(N_i)$. Since $N_i$ is incompressible,  $\tilde{N}_i$ is the
contractible universal cover of $N_i^{(m_i-1)}$.

By the previous Corollary, there are uniformly proper maps:
\[ \phi_1\colon\thinspace \tilde{N}_1 \vee \alpha \rw \tilde{W}_1 \]
\[ \phi_2\colon\thinspace \tilde{N}_2 \vee \beta \rw \tilde{W}_2 \]
So there exists a uniformly proper map:
\[ \phi_1\! \times \! \phi_2 \colon\thinspace (\tilde{N}_1 \vee \alpha) \! \times \! (\tilde{N}_2 \vee \beta) \rw
\tilde{W}_1 \!  \times \! \tilde{W}_2 \]
Recall that $(\tilde{N}_1 \vee \alpha) \! \times  \! (\tilde{N}_2 \vee \beta)$
is an $(m_1+m_2-1)$--proper obstructor by Proposition~\ref{join-lemma}.
Since $\tilde{W}_1  \! \times  \! \tilde{W}_2$ is quasi-isometric to $G_1  \! \times  \! G_2$\,:
\[ pobdim(G_1  \! \times  \! G_2) \geq m_1\! +\! m_2 \]
But $G_1  \! \times  \! G_2$ acts on $\tilde{W}_1 \! \times  \! \tilde{W}_2$ properly
discontinuously, and this implies that:
\[ pobdim(G_1 \!  \times  \! G_2) \leq actdim(G_1 \! \times  \! G_2) \leq m_1\! + \! m_2 \]
Therefore, $pobdim(G_1 \! \times  \! G_2)=m_1\! +\! m_2$.

Continue inductively and we conclude that:
\[ pobdim(G)=pobdim(G_1 \! \times   \cdots \times  \! G_d)= \! m_1\! +  \ldots  + \! m_d \]
Finally we see that
$$ pobdim(G) \leq updim(G) \leq actdim(G) \Rightarrow actdim(G)=\! m_1\! + \ldots + \! m_d \eqno{\qed}$$
\end{corollary}

\rk{Acknowledgements}
The author thanks Professor M. Bestvina for numerous helpful
discussions.


\appendix
\section{Pro-Category of Abelian Groups}
\label{pc}

With every category $\mathcal{K}$ we can associate a new category
$pro(\mathcal{K})$. We briefly review the definitions, see~\cite{AM}
or~\cite{Ma} for details.
Recall that a partially ordered set
$(\Lambda, \leq)$ is \emph{directed} if, for $i, j \in \Lambda$,
there exists $k \in \Lambda$ so that $k \geq i, j$.

\begin{definition}[(Inverse system)]
Let $(\Lambda, \leq)$ be a directed set.
The system $\mathbf{A}=\{A_{\lambda}, p^{\lambda'}_{\lambda}, \Lambda \}$ is
called an \emph{inverse system} over
$(\Lambda, \leq)$ in the category $\mathcal{K}$, if the following conditions
are true.
\begin{itemize}
\item[(i)] $A_{\lambda} \in Ob_{\mathcal{K}}$ for every $\lambda \in \Lambda$
\item[(ii)] $p^{\lambda'}_{\lambda} \in Mor_{\mathcal{K}}(A_{\lambda'}, A_{\lambda})$ for
$\lambda' \geq \lambda$
\item[(iii)] $\lambda \leq \lambda' \leq \lambda'' \Rightarrow
p^{\lambda'}_{\lambda}p^{\lambda''}_{\lambda'}=p^{\lambda''}_{\lambda}$
\end{itemize}
\end{definition}

\begin{definition}[(A map of systems)]
Given two inverse systems in $\mathcal{K}$,
\[\mathbf{A}=\{A_{\lambda}, p^{\lambda'}_{\lambda}, \Lambda \},\ \mbox{\ and\ }\
\mathbf{B}=\{B_{\mu}, q^{\mu'}_{\mu}, M \} \] the system
\[ \bar{f}=(f, f_{\lambda})\colon\thinspace\mathbf{A} \rightarrow \mathbf{B}\]
is called a \emph{map of systems} if the following conditions
are true.
\begin{itemize}
\item[(i)] $f\colon\thinspace M \rightarrow \Lambda$ is an increasing function
\item[(ii)] $f(M)$ is cofinal with $\Lambda$
\item[(iii)] $f_{\mu} \in Mor_{\mathcal{K}}(A_{f(\mu)}, B_{\mu})$
\item[(iv)] For $\mu' \geq \mu$ there exists $\lambda \geq f(\mu), f(\mu')$ so that:
\begin{gather*}
q^{\mu'}_{\mu} \circ f_{\mu} \circ p^{\lambda}_{f(\mu')} = f_{\mu} \circ p^{\lambda}_{\mu}\\
\begin{CD}
A_{f(\mu)} @<p^{f(\mu')}_{f(\mu)}<< A_{f(\mu')}  @<p^{\lambda}_{\mu}<<  A_{\lambda}\\
@Vf_{\mu}VV                        @Vf_{\mu'}VV \\
B_{\mu}  @<q^{\mu'}_{\mu}<<           B_{\mu'}
\end{CD}
\end{gather*}
\end{itemize}

\noindent
$\bar{f}$ is called a \emph{special map of systems}
if $\Lambda=M$, $f=id$, and
$f_{\lambda}p^{\lambda'}_{\lambda}=q^{\lambda'}_{\lambda}f_{\lambda'}$.
Two maps of systems $\bar{f}, \bar{g}\colon\thinspace\mathbf{A} \rightarrow \mathbf{B}$ are
considered equivalent, $\bar{f} \simeq \bar{g}$, if for every $\mu \in M$ there
is a $\lambda \in \Lambda, \lambda \geq f(\mu), g(\mu)$, such that
$f_{\mu}p^{\lambda}_{f(\mu)}=g_{\mu}p^{\lambda}_{g(\mu)}$.  This is an equivalence
relation.
\end{definition}

\begin{definition}[(Pro-category)]
$pro(\mathcal{K})$ is a category whose objects are inverse systems in $\mathcal{K}$
and morphisms are equivalence classes of maps of systems. The class
containing $\bar{f}$ will be denoted by $\mathbf{f}$.
\end{definition}

Our main interest is the following pro-category.

\begin{example}[Pro-category of abelian groups]
Let $\mathcal{A}$ be the category of abelian groups and
homomorphisms. Then corresponding $pro(\mathcal{A})$ is called the
category of pro-abelian groups.
\end{example}

\begin{example}[Homology pro-groups]
Suppose $\{{(X, X_0)}_i, p^{i'}_i, \nn \}$
is an object in the pro-homotopy category of pairs of spaces having
the homotopy type of a simplicial pair.
Then $\{H_j({(X, X_0)}_i), {(p^{i'}_{i})}_*, \nn \}$
is an object of $pro(\mathcal{A})$.  Denote
$\{H_j({(X, X_0)}_i), {(p^{i'}_{i})}_*, \nn \}$ by $proH_j(X, X_0)$.
\end{example}

We list useful properties of $pro(\mathcal{A})$:

\begin{enumerate}
\item A system $\mathbf{0}$ consisting of a single trivial group is a
zero object in $pro(\mathcal{A})$.
\item A pro-abelian group $\{G_i, p^{i'}_{i}, \nn \}
\cong \mathbf{0}$ iff every
$i$ admits a $i' \geq i$ such that $p^{i'}_{i}=0$.
\item Let $\mathbf{A}$ denote a constant pro-abelian group $\{A, id_{A}, \nn \}$.
If a pro-abelian group $\{G_i, p^{i'}_{i}, \nn \} \cong \mathbf{A}$
then \[\lim_{\lw} G_i = A. \]
See~\cite[Lemma 4.1]{EG}.
\end{enumerate}

\Addresses\recd
\end{document}